\documentclass[a4paper,11pt]{amsart}
\usepackage[lmargin=30mm,rmargin=30mm,tmargin=26mm,bmargin=26mm,footskip=7mm]{geometry}
\usepackage{amssymb}
\usepackage[dvipsnames]{xcolor}
\usepackage{microtype}
\usepackage[nocompress]{cite}
\usepackage[colorlinks=true,linkcolor=blue!50!black,citecolor=blue!50!black,urlcolor=blue!50!black,hyperfootnotes=false]{hyperref}
\usepackage{enumitem}
\usepackage{newunicodechar}
\newunicodechar{∀}{\ensuremath{\forall}}
\newunicodechar{ℕ}{\ensuremath{\mathbb{N}}}
\newunicodechar{ℝ}{\ensuremath{\mathbb{R}}}
\newunicodechar{≤}{\ensuremath{\leqslant}}
\newunicodechar{→}{\ensuremath{\to}}
\newunicodechar{√}{\ensuremath{\surd}}

\theoremstyle{plain}
\newtheorem{theorem}{Theorem}
\newtheorem{lemma}[theorem]{Lemma}

\makeatletter
\renewcommand{\paragraph}{\@startsection{paragraph}{4}{0pt}{1.5ex plus .5ex minus .2ex}{-1em}{\normalfont\bfseries}}
\makeatother
\DeclareMathOperator{\tr}{tr}
\newcommand{\doi}[1]{\href{https://doi.org/#1}{\nolinkurl{doi:#1}}}
\newcommand{\defn}[1]{\textcolor{Maroon}{\emph{#1}}}

\makeatletter
\def\@secnumfont{\bfseries}
\def\section{\@startsection{section}{1}\z@{.7\linespacing\@plus\linespacing}{.5\linespacing}{\normalfont\bfseries\centering}}
\def\@settitle{\begin{center}\baselineskip20\p@\relax\Large\bfseries\scshape\@title\end{center}}
\def\@setauthors{\begingroup\def\thanks{\protect\thanks@warning}\trivlist\centering\@topsep30\p@\relax\advance\@topsep by -\baselineskip\item\relax\author@andify\authors\def\\{\protect\linebreak}{\large\authors}\par\smallskip{\small School of Mathematics, Monash University, Melbourne, Australia\par\texttt{david.wood@monash.edu}}\endtrivlist\endgroup}
\makeatother
\begin{document}

\title{An explicit bound on the extremal function \linebreak for complete graph minors}
\author{David R. Wood}
\thanks{Research supported by the Australian Research Council and by NSERC.}
\date{30 October 2026}

\begin{abstract}
Kostochka and Thomason independently proved that every graph with density at least $\beta\,t\sqrt{\log t}$ has a $K_t$ minor, for some constant $\beta$. Thomason later proved that the least such density is $(0.319\ldots+o(1))\,t\sqrt{\log t}$. However, the $o(1)$ term here is not explicit, so the formula gives no bound for any specific value of $t$. The goal of this paper is a bound of the form $\beta\,t\sqrt{\log t}$ that holds for every $t$, with $\beta$ small. We prove that, for every integer $t\geqslant2$, every graph with density at least $4.5\,t\sqrt{\log t}$ has a $K_t$ minor.
\end{abstract}

\maketitle
\makeatletter
\def\@oddhead{}\def\@evenhead{}
\def\@oddfoot{\normalfont\hfil\thepage\hfil}\def\@evenfoot{\@oddfoot}
\makeatother

%=====================================================================
\section{Introduction}\label{sec:intro}
%=====================================================================

Mader~\cite{Mader67} first proved that every sufficiently dense graph contains the complete graph $K_t$ as a minor. Here the \defn{density} of a graph $G$ is $d(G):=e(G)/v(G)$ if $v(G)\geqslant1$, and $d(G):=0$ if $v(G)=0$.\footnote{Graphs are finite and simple. The \defn{null} graph is the graph with no vertices. A graph $H$ is a \defn{minor} of a graph $G$ if a graph isomorphic to $H$ can be obtained from a subgraph of $G$ by contracting edges. For a graph $G$, let $v(G)$ and $e(G)$ be the numbers of vertices and edges of $G$, respectively. We write $\log$ for the natural logarithm, $\log_2$ for the logarithm to base $2$, and $\mathrm{e}$ for Euler's number.} For an integer $t\geqslant1$, let $c(K_t)$ be the least real number $c$ such that every graph $G$ with $d(G)>c$ has a $K_t$ minor. Kostochka~\cite{Kostochka82,Kostochka84} and Thomason~\cite{Thomason84} independently proved that $c(K_t)\in O(t\sqrt{\log t})$, which is tight for suitable random graphs~\cite{BCE80,FdlV83,Kostochka84}. Thomason~\cite{Thomason01} later proved that $c(K_t)=(\alpha+o(1))\,t\sqrt{\log t}$, where $\alpha=0.3190863\ldots$. However, the $o(1)$ term in Thomason's formula is not explicit, so it gives no upper bound on $c(K_t)$ for specific values of $t$.

At the \emph{Probabilistic and Extremal Combinatorics Downunder} workshop, held at Monash University in 2016, I had the pleasure of discussing the extremal function for excluded complete graph minors with Andrew Thomason and Alexandr Kostochka, who both said that it would be desirable to determine a small value of $\beta$ such that, for every integer $t$, including small values, every $K_t$-minor-free graph has density at most $\beta\,t\sqrt{\log t}$. The probabilistic method handles large values of $t$ well, but medium values of $t$ require significant computation, so I employed Claude Opus 5.5 to help out, and this collaboration led to the following result.

\begin{theorem}\label{thm:Kt}
For every integer $t\geqslant2$, every graph $G$ with $d(G)\geqslant4.5\,t\sqrt{\log t}$ has a $K_t$ minor, and so $c(K_t)\leqslant4.5\,t\sqrt{\log t}$.
\end{theorem}

\paragraph{Statement of AI Use.} I started by asking Claude Opus 5.5 to optimise the constant in the proof of Thomason~\cite{Thomason84}. Claude identified several areas for optimisation, worked out the details, and produced a first draft of the paper. I thoroughly revised the whole paper to improve the exposition. I verified the proof at a high level, but have not checked all the technical and numerical details by hand. For these details I rely on automated verification. Claude wrote the numeric scripts and the Lean formalisation.

\paragraph{Related work.} For small $t$, the exact value of $c(K_t)$ is known. Vizing~\cite{Vizing68} asked for the maximum number of edges in a graph with $n$ vertices and no $K_t$ minor. Zelinka~\cite{Zelinka73} conjectured that this maximum is $(t-2)n-\binom{t-1}2$ for $n\geqslant t-1$. Dirac~\cite{Dirac64} and Mader~\cite{Mader68} respectively proved this conjecture for $t\leqslant5$ and $t\leqslant7$, and Mader observed that the conjecture fails for $t=8$, and hence for all $t\geqslant8$. Zelinka~\cite{Zelinka73} and Gy\H{o}ri~\cite{Gyori82} respectively gave independent proofs for $t\leqslant4$ and for $t\in\{5,6\}$. J{\o}rgensen~\cite{Jorgensen94}, Song and Thomas~\cite{SongThomas06} and Zhu~\cite{Zhu21} respectively determined the extremal function for $K_8$, $K_9$ and $K_{10}$ minors. These results imply that $c(K_t)\leqslant t-2$ for $t\in\{2,\dots,10\}$. In fact, $c(K_t)=t-2$ for $t\in\{2,\dots,10\}$, since the graph obtained from $K_{t-2}$ by adding $n$ pairwise non-adjacent vertices, each adjacent to every vertex of $K_{t-2}$, has no $K_t$ minor, and the density of this graph tends to $t-2$ as $n\to\infty$.

The first general upper bounds on $c(K_t)$ are due to Mader~\cite{Mader67}, who proved that $c(K_t)\leqslant2^{t-3}$ for every $t\geqslant3$. Mader~\cite{Mader68} later improved this bound to $c(K_t)\leqslant8(t-2)\lfloor\log_2(t-2)\rfloor$ for $t\geqslant4$. As mentioned above, Kostochka~\cite{Kostochka82,Kostochka84} and Thomason~\cite{Thomason84,Thomason01} improved the bound further, to $O(t\sqrt{\log t})$, and this order of growth is best possible~\cite{BCE80,FdlV83,Kostochka84}. Alon, Krivelevich and Sudakov~\cite{AKS23} recently gave a new proof that $c(K_t)\in O(t\sqrt{\log t})$. Delcourt and Postle~\cite[Theorem~4.3]{DP25} state that it follows from the work of Kostochka~\cite{Kostochka82} that for every integer $t\geqslant2$, every graph $G$ with $d(G)\geqslant30\,t\sqrt{\log t}$ has a $K_t$ minor. To my knowledge, this bound with constant $30$ is the best explicit bound of the form $c(K_t)\leqslant\beta\,t\sqrt{\log t}$ that holds for every $t\geqslant2$. Theorem~\ref{thm:Kt} improves the constant $30$ to $4.5$.

Extremal functions for $H$-minor-free graphs have been widely studied (see the survey~\cite{Thomason06}), including when $H$ is a graph with $t$ vertices and average degree $d$~\cite{MT05,RW16,HarveyWood16,NRTW20,TW22,HNW22,HKL22}, a disconnected graph such as a union of cycles~\cite{Justesen88,HarveyWood15,CLNWY17,KNQ21,BohmeKostochka05,Thomason07,Thomason08}, a complete bipartite graph~\cite{Myers03,KO05,KP08,KP10,KP12,CRS11,Wales21}, a complete graph minus one or two edges~\cite{Dirac64,Jakobsen71,Jakobsen72,Jakobsen83,Song04,Song05,Rolek20}, a complete graph~\cite{CDGV07}, a path or a cycle~\cite{ErdosGallai59}, $K_{4,4}$~\cite{Jorgensen01}, $K_{3,3}$~\cite{Hall43}, the octahedron $K_{2,2,2}$~\cite{Ding13}, or the Petersen graph~\cite{HendreyWood18}. Other work restricts the host graph to be triangle-free, bipartite, apex or of given girth~\cite{ThomasYoo19,Pavelescu22,McCartyThomas19,CSSS24,Barat23,Barat24}. Related work studies the possible densities of minor-closed classes~\cite{Eppstein10,MP19,KapadiaNorin20}. Hendrey, Norin, Steiner and Turcotte~\cite{HNST25} studied how dense a $t$-vertex minor must be in every graph with average degree at least $t-1$.

\paragraph{Applications.} The maximum density of $K_t$-minor-free graphs is of fundamental importance in structural and algorithmic graph theory. For example, a minimum-degree greedy algorithm shows that every $K_t$-minor-free graph is $(\lfloor2c(K_t)\rfloor+1)$-colourable, and is in fact $(\lfloor2c(K_t)\rfloor+1)$-choosable. So, by Theorem~\ref{thm:Kt}, every $K_t$-minor-free graph is $\lceil9\,t\sqrt{\log t}\,\rceil$-choosable. Hadwiger's Conjecture~\cite{Hadwiger43} asserts that every $K_t$-minor-free graph is $(t-1)$-colourable. For several decades, the best known upper bound on the chromatic number of $K_t$-minor-free graphs was $O(t\sqrt{\log t})$, which follows from the Kostochka--Thomason bound on $c(K_t)$. Norin, Postle and Song~\cite{NPS23} were the first to improve this order of magnitude, Delcourt and Postle~\cite{DP25} improved the bound to $O(t\log\log t)$, and Norin and Steiner~\cite{NorinSteiner26} recently proved that $O(t)$ colours suffice (see also~\cite{Wood26,OpenAI26}). The Kostochka--Thomason bound remains a key ingredient in several of these recent papers, for example in~\cite[Theorem~4.3]{DP25}, so Theorem~\ref{thm:Kt} can be used to improve the constants. As another application, Norine, Seymour, Thomas and Wollan~\cite{NSTW06} used the Kostochka--Thomason bound to show that every proper minor-closed class contains at most $c^n\,n!$ graphs with vertex set $\{1,\dots,n\}$, for some constant $c$ depending on the class. The Kostochka--Thomason bound also has algorithmic applications, via arboricity. By a theorem of Nash-Williams~\cite{NashWilliams64}, the Kostochka--Thomason bound implies that every $K_t$-minor-free graph has arboricity $O(t\sqrt{\log t})$, that is, its edge set can be partitioned into $O(t\sqrt{\log t})$ forests. Many distributed and parallel algorithms for colouring, maximal independent sets and dominating sets have running times or approximation ratios that depend on the arboricity of the input graph~\cite{GPS88,BarenboimElkin10,LenzenWattenhofer10}, so Theorem~\ref{thm:Kt} makes this dependence explicit for $K_t$-minor-free graphs.

\paragraph{Proof Sketch.} The proof of Theorem~\ref{thm:Kt} follows the argument of Thomason~\cite{Thomason84} with explicit constants. Two tools come from that paper. The first, Lemma~\ref{lem:red2}, passes to a minor that is minimal subject to a lower bound on its number of edges. The second is the core of the proof (Lemma~\ref{lem:core}), which builds each branch set from a few random vertices and shows, by Markov's inequality and the union bound, that the branch sets can be chosen to be pairwise adjacent. The proof also uses two ideas from the later paper of Thomason~\cite{Thomason01}. The first idea is to exploit the actual edge density of the dense minor produced by the reductions, instead of only the fact that every vertex of this minor is adjacent to more than half of the other vertices. In this paper, this idea leads to the density functions in the second reduction and to the bounds $\pi_i$ on the proportion of non-neighbours in the core of the proof. The second idea is to connect a set of vertices through a vertex that dominates most of the set; see~\cite[Lemma~4.1]{Thomason01} and Lemma~\ref{lem:connect} in this paper.

Let $G$ be a graph with $d(G)\geqslant c$, where $c\leqslant4.5\,t\sqrt{\log t}$ is chosen for each $t$. Section~\ref{sec:reductions} applies minor-minimality twice. First, as in Mader's proof that $c(K_t)\leqslant2^{t-3}$~\cite{Mader67}, Lemma~\ref{lem:red1} takes a minor $G_1$ of $G$ with $d(G_1)\geqslant c$ and $v(G_1)+e(G_1)$ minimum, and lets $H_2$ be the subgraph of $G_1$ induced by the neighbourhood of a vertex of minimum degree in $G_1$. Second, following Thomason~\cite{Thomason84}, Lemma~\ref{lem:red2} takes a minor $J$ of $H_2$ with $v(J)+e(J)$ minimum subject to $e(J)\geqslant\Phi(v(J))$, where $\Phi$ is a suitable density function, and lets $G_0$ be the subgraph of $J$ induced by the neighbourhood of a vertex of minimum degree in $J$. The graph $G_0$ has bounded order, and $2\delta(G_0)-v(G_0)$ is large, with bounds that depend on $v(J)$. Each of the two reductions finds one vertex of the complete minor (Lemma~\ref{lem:apex}), so it suffices to find a $K_{t-2}$ minor in $G_0$.

Section~\ref{sec:core} builds the branch sets of a $K_{t-2}$ minor in $G_0$ one at a time (Lemma~\ref{lem:core}). Since $2\delta(G_0)-v(G_0)$ is large, every vertex of $G_0$ is adjacent to more than half of the other vertices, and this remains true after the earlier branch sets are removed. More precisely, the bounds on $2\delta(G_0)-v(G_0)$ and $v(G_0)$ give an upper bound $\pi_i<\frac12$ on the proportion of non-neighbours of each vertex when the $i$-th branch set is chosen. So a few random vertices are adjacent to almost all of the remaining vertices. For $i\in\{1,\dots,t-3\}$, the $i$-th branch set consists of $k_i$ random vertices and a few further vertices that make the branch set connected. Each new branch set must be adjacent to every earlier branch set, and this holds with positive probability, since each earlier branch set is adjacent to almost all of the remaining vertices. The last branch set consists of all remaining vertices. Larger values of $k_i$ make the $i$-th branch set adjacent to more vertices, which helps the later branch sets, but larger branch sets use up more of the vertices of $G_0$, and the choice of the sequence $(k_i)$ balances these two effects. The bound $\pi_i$ increases as the construction proceeds, so the early branch sets can start with fewer random vertices than the later ones. Moreover, at most $\lfloor k_i\pi_i\rfloor+1$ further vertices make the $i$-th branch set connected (Lemma~\ref{lem:connect}).

Section~\ref{sec:proof} combines the lemmas of Sections~\ref{sec:reductions} and~\ref{sec:core} (Lemma~\ref{lem:assembly}) and handles three ranges of $t$. For $t\in\{4,\dots,400\}$, a computer finds a density function and, for each possible value of $v(J)$, a sequence $(k_i)$ (Lemma~\ref{lem:small}). For $t\in\{401,\dots,10^{35}\}$, the density function is given by a formula with several slopes, and $k_i$ is determined by the interval that contains $\pi_i$. A computer checks the resulting inequalities for finitely many cells, each of which covers an interval of values of $t$ (Lemma~\ref{lem:middle}). For $\log t\geqslant80$, an analytic argument handles every such $t$ (Lemma~\ref{lem:large}). So the proof of Theorem~\ref{thm:Kt} relies on computer checks in exact and interval arithmetic, which are described in Appendix~\ref{sec:verification}.

\paragraph{Notation.} We use standard graph-theoretic notation. For a graph $G$ and a vertex $u$ of $G$, let $N_G(u)$ be the neighbourhood of $u$, let $\deg_G(u):=|N_G(u)|$, and, if $G$ is a non-null graph, let $\delta(G)$ be the minimum degree of $G$. For $Y\subseteq V(G)$, let $G[Y]$ be the subgraph of $G$ induced by $Y$, and let $G-Y:=G[V(G)\setminus Y]$. The \defn{closed neighbourhood} $N_G[Y]$ is the set of vertices of $G$ that are in $Y$ or adjacent to a vertex of $Y$, and $N_G[u]:=N_G[\{u\}]$. A \defn{non-neighbour} of $u$ is a vertex $w\neq u$ that is not adjacent to $u$. For an edge $uw$ of $G$, let $\tr_G(uw):=|N_G(u)\cap N_G(w)|$ be the number of triangles of $G$ that contain $uw$. We drop the subscript $G$ when $G$ is clear from the context. For every integer $n$, let $\binom n2:=n(n-1)/2$.

For an integer $s\geqslant1$, a \defn{model} of $K_s$ in a graph $G$ is a family $(X_1,\dots,X_s)$ of pairwise disjoint non-empty sets $X_i\subseteq V(G)$, called \defn{branch sets}, such that each $G[X_i]$ is connected and, for all distinct $i$ and $j$, some edge of $G$ joins $X_i$ and $X_j$. A graph has a $K_s$ minor if and only if it has a model of $K_s$.

%=====================================================================
\section{Two reductions}\label{sec:reductions}
%=====================================================================

Lemma~\ref{lem:red1} goes back to Mader~\cite{Mader67}, and Lemma~\ref{lem:red2} is based on Lemma~1 of Thomason~\cite{Thomason84}. We include proofs of both lemmas to keep track of integer rounding.

\begin{lemma}\label{lem:apex}
For every integer $s\geqslant1$ and every vertex $u$ of a graph $J$, if $J[N_J(u)]$ has a $K_s$ minor, then $J$ has a $K_{s+1}$ minor.
\end{lemma}

\begin{proof}
If $(X_1,\dots,X_s)$ is a model of $K_s$ in $J[N_J(u)]$, then $(X_1,\dots,X_s,\{u\})$ is a model of $K_{s+1}$ in $J$.
\end{proof}

\begin{lemma}\label{lem:red1}
For every real number $c>0$, every graph $G$ with $d(G)\geqslant c$ has a minor $G_1$ with a vertex $z$ such that the graph $H_2:=G_1[N_{G_1}(z)]$ satisfies
\[
c<v(H_2)<2c+1\qquad\text{and}\qquad\delta(H_2)\geqslant\lfloor c\rfloor ,
\]
and such that, for every integer $s\geqslant1$, if $H_2$ has a $K_s$ minor, then $G$ has a $K_{s+1}$ minor.
\end{lemma}

\begin{proof}
Let $G_1$ be a minor of $G$ with $d(G_1)\geqslant c$ and $v(G_1)+e(G_1)$ minimum, and let $n_1:=v(G_1)$. Since $d(G_1)\geqslant c>0$, the graph $G_1$ has an edge, so $n_1\geqslant2$. By minimality, every minor of $G_1$ obtained by deleting a vertex, deleting an edge or contracting an edge has density less than $c$. In particular, for an edge $uw$, contracting $uw$ gives $e(G_1)-1-\tr(uw)<c(n_1-1)$. Since $e(G_1)\geqslant c\,n_1$, it follows that $\tr(uw)>c-1$. Since $\tr(uw)$ is an integer, $\tr(uw)\geqslant\lfloor c\rfloor$. Similarly, deleting an edge gives $e(G_1)-1<c\,n_1$.

Let $z$ be a vertex of minimum degree in $G_1$. Then $v(H_2)=\deg(z)\leqslant2e(G_1)/n_1<2c+2/n_1\leqslant2c+1$. Deleting $z$ from $G_1$ gives $e(G_1)-\deg(z)<c(n_1-1)$, and so $\deg(z)>c$. Each vertex $w\in N(z)$ has $\tr(zw)\geqslant\lfloor c\rfloor$ neighbours in $H_2$, so $\delta(H_2)\geqslant\lfloor c\rfloor$. The final claim follows from Lemma~\ref{lem:apex}, since $G_1$ is a minor of $G$.
\end{proof}

The second reduction uses the following functions. Let $c>0$ be a real number, and let $N:=\lceil2c+1\rceil-1$, which is the largest integer less than $2c+1$. Since $\lfloor c\rfloor<2c$, we have $N\geqslant\lfloor c\rfloor+1$. Let $n_0\leqslant\lfloor c\rfloor$ be a positive integer. A \defn{density function for $c$ and $n_0$} is a non-decreasing function $\Phi:\{n_0,\dots,N\}\to\mathbb Z$ such that $\Phi(n_0)>\binom{n_0}2$ and
\begin{equation}\label{eq:density}
\Phi(n)\ \leqslant\ \Bigl\lceil\frac{n\lfloor c\rfloor}2\Bigr\rceil\qquad\text{for every }n\in\{\lfloor c\rfloor+1,\dots,N\} .
\end{equation}
For such a function and every $n\in\{n_0+1,\dots,N\}$, let
\[
V(n)\ :=\ \Bigl\lfloor\frac{2\Phi(n)}n\Bigr\rfloor\qquad\text{and}\qquad D(n)\ :=\ \Phi(n)-\Phi(n-1).
\]
The number $V(n)$ is an upper bound on the minimum degree of every graph with $n$ vertices and $\Phi(n)$ edges, and so $V(n)$ bounds the number of vertices of the graph $G_0$ in Lemma~\ref{lem:red2} below.

\begin{lemma}\label{lem:red2}
For every real number $c>0$, every positive integer $n_0\leqslant\lfloor c\rfloor$ and every density function $\Phi$ for $c$ and $n_0$, every graph $H_2$ with $c<v(H_2)<2c+1$ and $\delta(H_2)\geqslant\lfloor c\rfloor$ has a minor $J$ with a vertex $u$ such that $n:=v(J)$ satisfies $n\in\{n_0+1,\dots,N\}$ and $\Phi(n)\leqslant\binom n2$, the graph $G_0:=J[N_J(u)]$ is a non-null graph with $v(G_0)\leqslant V(n)$ and $\delta(G_0)\geqslant D(n)$, and, for every integer $s\geqslant1$, if $G_0$ has a $K_s$ minor, then $H_2$ has a $K_{s+1}$ minor.
\end{lemma}

\begin{proof}
Let $\mathcal D$ be the set of graphs $L$ with $n_0\leqslant v(L)\leqslant N$ and $e(L)\geqslant\Phi(v(L))$. Let $n_2:=v(H_2)$. Since $n_2$ is an integer with $c<n_2<2c+1$, we have $n_2\in\{\lfloor c\rfloor+1,\dots,N\}$. Since $\delta(H_2)\geqslant\lfloor c\rfloor$ and $e(H_2)$ is an integer, $e(H_2)\geqslant\lceil n_2\lfloor c\rfloor/2\rceil\geqslant\Phi(n_2)$ by~\eqref{eq:density}. So $H_2\in\mathcal D$.

Let $J$ be a minor of $H_2$ in $\mathcal D$ with $v(J)+e(J)$ minimum, and let $n:=v(J)$. Then $n_0\leqslant n\leqslant n_2\leqslant N$. If $n=n_0$, then $e(J)\leqslant\binom{n_0}2<\Phi(n_0)$, which contradicts $J\in\mathcal D$. So $n\geqslant n_0+1$. Every graph $J'$ obtained from $J$ by deleting a vertex or contracting an edge is a minor of $H_2$ with $n_0\leqslant v(J')=n-1\leqslant N$, and so $e(J')<\Phi(n-1)$ by minimality. Since $\Phi$ is non-decreasing, $\Phi(n)\geqslant\Phi(n_0)\geqslant1$, and so $J$ has an edge. Deleting an edge of $J$ gives $e(J)-1<\Phi(n)$ by minimality. Since $e(J)\geqslant\Phi(n)$, it follows that $e(J)=\Phi(n)$, and in particular $\Phi(n)\leqslant\binom n2$. In particular, deleting a vertex $x$ of $J$ gives $e(J)-\deg(x)<\Phi(n-1)$, and so $\deg(x)>D(n)$. Similarly, contracting an edge $xw$ of $J$ gives $e(J)-1-\tr(xw)<\Phi(n-1)$, and so $\tr(xw)>D(n)-1$. Since $\tr(xw)$ and $D(n)$ are integers, $\tr(xw)\geqslant D(n)$.

Let $u$ be a vertex of minimum degree in $J$. Then $v(G_0)=\deg(u)>D(n)\geqslant0$, where $D(n)\geqslant0$ since $\Phi$ is non-decreasing. So $G_0$ is a non-null graph. Also $v(G_0)\leqslant2e(J)/n=2\Phi(n)/n$, and since $v(G_0)$ is an integer, $v(G_0)\leqslant V(n)$. Each vertex $w$ of $G_0$ has $\tr(uw)\geqslant D(n)$ neighbours in $G_0$, so $\delta(G_0)\geqslant D(n)$. The final claim follows from Lemma~\ref{lem:apex}, since $J$ is a minor of $H_2$.
\end{proof}

%=====================================================================
\section{Core lemma}\label{sec:core}
%=====================================================================

This section proves the core lemma (Lemma~\ref{lem:core}), which finds a $K_m$ minor in every non-null graph $G_0$ for which $2\delta(G_0)-v(G_0)$ is sufficiently large and $v(G_0)$ is bounded, such as the graph $G_0$ given by Lemma~\ref{lem:red2}. We first give two properties of graphs in which every vertex has few non-neighbours.

\begin{lemma}\label{lem:dense}
For every integer $q\geqslant0$ and every graph $G$ with $n\geqslant2q+1$ vertices in which every vertex has at most $q$ non-neighbours, any two distinct non-adjacent vertices of $G$ have a common neighbour, and so $G$ is connected.
\end{lemma}

\begin{proof}
Two distinct non-adjacent vertices each have at least $n-1-q$ neighbours among the other $n-2$ vertices, so these two vertices have at least $2(n-1-q)-(n-2)=n-2q\geqslant1$ common neighbours. So $G$ has diameter at most $2$, and in particular $G$ is connected.
\end{proof}

The branch sets in the core lemma have sizes bounded by the following numbers. For an integer $k\geqslant1$ and a real number $\pi\geqslant0$, let
\[
\ell_1(\pi)\ :=\ 1\qquad\text{and}\qquad\ell_k(\pi)\ :=\ k+1+\lfloor k\pi\rfloor\quad\text{for }k\geqslant2 .
\]

\begin{lemma}\label{lem:connect}
For all $q$, $G$ and $n$ as in Lemma~\ref{lem:dense}, every real number $\pi\geqslant q/n$, every integer $k\geqslant1$ and every non-empty set $Y\subseteq V(G)$ with $|Y|\leqslant k$, there is a set $X$ with $Y\subseteq X\subseteq V(G)$ such that $G[X]$ is connected and $|X|\leqslant\ell_k(\pi)$.
\end{lemma}

\begin{proof}
Let $j:=|Y|$. If $j=1$, then $X:=Y$ satisfies the lemma, since $\pi\geqslant q/n\geqslant0$ and so $\ell_k(\pi)\geqslant1$. Now let $j\geqslant2$, so that $k\geqslant2$. Since every vertex of $G$ has degree at least $n-1-q$,
\[
\sum_{z\in V(G)}|N[z]\cap Y|\ =\ \sum_{y\in Y}\bigl(\deg(y)+1\bigr)\ \geqslant\ j\,(n-q).
\]
So some vertex $z$ of $G$ satisfies $|N[z]\cap Y|\geqslant j(n-q)/n$. Then $Q:=Y\setminus N[z]$ satisfies $|Q|\leqslant jq/n\leqslant k\pi$. Since $|Q|$ is an integer, $|Q|\leqslant\lfloor k\pi\rfloor$. By Lemma~\ref{lem:dense}, each $y\in Q$ has a common neighbour $x_y$ with $z$. Then $X:=Y\cup\{z\}\cup\{x_y:y\in Q\}$ induces a connected subgraph of $G$, since every vertex of $X$ is $z$, a neighbour of $z$, or a vertex $y\in Q$ adjacent to the neighbour $x_y$ of $z$. Finally, $|X|\leqslant j+1+|Q|\leqslant k+1+\lfloor k\pi\rfloor=\ell_k(\pi)$.
\end{proof}

The hypotheses of the core lemma use the following numbers. Let $m\geqslant2$ and $R<V$ be integers, and let $k_1,\dots,k_{m-1}$ and $L_1,\dots,L_{m-1}$ be positive integers with $\sum_{i=1}^{m-1}L_i\leqslant R$. We call $L_i$ the \defn{charged cost} of position~$i$ of the sequence $(k_i)$. For $j\in\{0,\dots,m-1\}$, let $S_j:=\sum_{a=1}^{j}L_a$. For $i\in\{1,\dots,m\}$ and $a\in\{1,\dots,m-1\}$, let
\[
\pi_i:=\frac{V-R-1}{2(V-S_{i-1})}\qquad\text{and}\qquad w_a:=\Bigl(1+\frac1{k_a}\Bigr)\pi_a^{\,k_a-1}.
\]
Since $0\leqslant S_{i-1}\leqslant R<V$, we have $0\leqslant\pi_1\leqslant\pi_2\leqslant\dots\leqslant\pi_m<\frac12$. In the proof of the core lemma, $\pi_i$ bounds the proportion of non-neighbours of each vertex when the $i$-th branch set is chosen. For an integer $k\geqslant1$, let $\psi(k):=(1+\frac1k)2^{-k}$, and note that $\psi(k)\leqslant1$. When the sequence $(k_i)$ is constant, the proof of the core lemma is the counting argument of Thomason~\cite{Thomason84}, with explicit constants.

\begin{lemma}\label{lem:core}
For all integers $m\geqslant2$ and $R<V$, and all positive integers $k_1,\dots,k_{m-1}$ and $L_1,\dots,L_{m-1}$ with $\sum_{i=1}^{m-1}L_i\leqslant R$ such that $L_i\geqslant\ell_{k_i}(\pi_i)$ for $i\in\{1,\dots,m-1\}$ and such that
\begin{equation}\label{eq:core}
\sum_{a=1}^{i-1}\bigl(w_a\,\pi_i\bigr)^{k_i}\ <\ \frac1{k_i+1}\qquad\text{for }i\in\{1,\dots,m-1\} ,
\end{equation}
every non-null graph $G_0$ with $v(G_0)\leqslant V$ and $2\delta(G_0)-v(G_0)\geqslant R-1$ has a $K_m$ minor.
\end{lemma}

\begin{proof}
Let $n:=v(G_0)$ and $q:=n-1-\delta(G_0)$. Then every vertex of $G_0$ has at most $q$ non-neighbours, and the hypothesis $2\delta(G_0)-n\geqslant R-1$ says that $2q+1\leqslant n-R$. Since $L_i\geqslant1$ for every $i$, we have $R\geqslant m-1$, and so $n\geqslant R+1\geqslant m$. If $q=0$, then $G_0$ is a complete graph with at least $m$ vertices, and so $G_0$ has a $K_m$ minor. So we may assume that $q\geqslant1$.

We construct pairwise disjoint non-empty sets $X_1,\dots,X_m\subseteq V(G_0)$ such that each $G_0[X_i]$ is connected, $|X_i|\leqslant L_i$ for $i\leqslant m-1$, and some edge of $G_0$ joins $X_a$ and $X_i$ for all $a<i$. Then $(X_1,\dots,X_m)$ is a model of $K_m$.

Let $j\in\{0,\dots,m-1\}$, and suppose that $X_1,\dots,X_j$ have been constructed. Let $G':=G_0-(X_1\cup\dots\cup X_j)$ and $n':=v(G')$. Then $n'\geqslant n-S_j\geqslant n-R\geqslant2q+1$, and every vertex of $G'$ has at most $q$ non-neighbours in $G'$. So Lemmas~\ref{lem:dense} and~\ref{lem:connect} apply to $G'$, with $n'$ in place of $n$. Since $y\mapsto(y-R-1)/(y-S_j)$ is increasing on $(S_j,\infty)$, $2q\leqslant n-R-1$ and $n\leqslant V$,
\begin{equation}\label{eq:half}
\frac q{n'}\ \leqslant\ \frac{n-R-1}{2(n-S_j)}\ \leqslant\ \frac{V-R-1}{2(V-S_j)}\ =\ \pi_{j+1}\ <\ \frac12 .
\end{equation}
For $a\leqslant j$, let $G^{(a)}:=G_0-(X_1\cup\dots\cup X_{a-1})$, which is the graph in which $X_a$ was chosen, let $n^{(a)}:=v(G^{(a)})$, and let $U_a:=V(G^{(a)})\setminus N_{G^{(a)}}[X_a]$. We maintain the invariant
\begin{equation}\label{eq:inv}
|U_a|\ \leqslant\ \Bigl(1+\frac1{k_a}\Bigr)q^{k_a}\bigl(n^{(a)}\bigr)^{1-k_a}\qquad\text{for every }a\leqslant j.
\end{equation}
Let $a\leqslant j$. Inequality~\eqref{eq:half} at the step at which $X_a$ was chosen gives $q/n^{(a)}\leqslant\pi_a\leqslant\frac12$. Combining this bound with~\eqref{eq:inv}, and with~\eqref{eq:half} for the current value of $j$, gives
\begin{equation}\label{eq:Ua}
\frac{|U_a\cap V(G')|}{n'}\ \leqslant\ \Bigl(1+\frac1{k_a}\Bigr)\Bigl(\frac q{n^{(a)}}\Bigr)^{k_a-1}\frac q{n'}\ \leqslant\ w_a\,\pi_{j+1}\ <\ \Bigl(1+\frac1{k_a}\Bigr)2^{-k_a}\ =\ \psi(k_a).
\end{equation}
Every vertex of $V(G')\setminus U_a$ lies in $N_{G^{(a)}}[X_a]\setminus X_a$, and so is adjacent to a vertex of $X_a$.

First suppose that $j=m-1$, and let $X_m:=V(G')$. Then $G'[X_m]$ is connected by Lemma~\ref{lem:dense}. For each $a<m$, inequality~\eqref{eq:Ua} and $\psi(k_a)\leqslant1$ give $|U_a\cap V(G')|<n'$. So some vertex of $X_m$ is not in $U_a$, and this vertex is adjacent to a vertex of $X_a$.

Now suppose that $j\leqslant m-2$. We construct the branch set $X_i$, where $i:=j+1$. Let $k:=k_i$. Choose vertices $y_1,\dots,y_k\in V(G')$ independently and uniformly at random, and let $Y_i:=\{y_1,\dots,y_k\}$. Each vertex $w$ of $G'$ has at most $q$ non-neighbours, so $w\notin N_{G'}[Y_i]$ with probability at most $(q/n')^k$. Since $q\geqslant1$, Markov's inequality shows that the event
\[
\bigl|V(G')\setminus N_{G'}[Y_i]\bigr|\ >\ \Bigl(1+\frac1k\Bigr)q^k\,(n')^{1-k}
\]
has probability at most $\frac k{k+1}$. For each $a<i$, by~\eqref{eq:Ua}, the event $Y_i\subseteq U_a$ has probability at most $(w_a\pi_i)^k$. By~\eqref{eq:core} and the union bound, with positive probability none of these $i$ events occurs. Fix an outcome $Y_i$ for which none of these events occurs. By~\eqref{eq:half}, $q/n'\leqslant\pi_i$. So Lemma~\ref{lem:connect} with $\pi:=\pi_i$, applied to $G'$ and $Y_i$, gives a set $X_i\supseteq Y_i$ such that $G'[X_i]$ is connected and $|X_i|\leqslant\ell_k(\pi_i)\leqslant L_i$. Since $X_i\subseteq V(G')$, the set $X_i$ is disjoint from $X_1\cup\dots\cup X_j$, and $G_0[X_i]=G'[X_i]$ is connected. Since $G^{(i)}=G'$ and $Y_i\subseteq X_i$, we have $U_i\subseteq V(G')\setminus N_{G'}[Y_i]$. Since the displayed event does not occur, \eqref{eq:inv} holds for $a=i$. For each $a<i$, since $Y_i\not\subseteq U_a$, some vertex of $Y_i$ is not in $U_a$, and this vertex is adjacent to a vertex of $X_a$.
\end{proof}

Given $k_1,\dots,k_{m-1}$, $R$ and $V$, the \defn{canonical costs} are defined recursively by $L_i:=\ell_{k_i}(\pi_i)$, where $\pi_i$ is computed from $L_1,\dots,L_{i-1}$. The canonical cost $L_i$ is defined as long as $S_{i-1}<V$.

%=====================================================================
\section{\texorpdfstring{Proof of Theorem~\ref{thm:Kt}}{Proof of Theorem 1}}\label{sec:proof}
%=====================================================================

This section proves Theorem~\ref{thm:Kt} by combining the reductions of Section~\ref{sec:reductions} with the core lemma (Lemma~\ref{lem:assembly}).

Let $t\geqslant4$ be an integer. For integers $V$ and $R$, a sequence $k_1,\dots,k_{t-3}$ of positive integers with charged costs $L_1,\dots,L_{t-3}$ is \defn{admissible for $t$ with respect to $V$ and $R$} if the sequence satisfies the hypotheses of Lemma~\ref{lem:core} with $m=t-2$, that is, $\sum_{i=1}^{t-3}L_i\leqslant R<V$, $L_i\geqslant\ell_{k_i}(\pi_i)$ for every $i$, and~\eqref{eq:core} holds. A pair $(V,D)$ of integers is \defn{certified for $t$} if $D\geqslant V$, or if there is an integer $D'\leqslant D$ and a sequence with charged costs that is admissible for $t$ with respect to $V$ and $2D'-V+1$.

\begin{lemma}\label{lem:assembly}
For every integer $t\geqslant4$, every real number $c>0$, every positive integer $n_0\leqslant\lfloor c\rfloor$ and every density function $\Phi$ for $c$ and $n_0$ such that $(V(n),D(n))$ is certified for $t$ for every $n\in\{n_0+1,\dots,N\}$ with $\Phi(n)\leqslant\binom n2$, every graph $G$ with $d(G)\geqslant c$ has a $K_t$ minor.
\end{lemma}

\begin{proof}
By Lemma~\ref{lem:red1}, $G$ has a minor $G_1$ with a vertex $z$ such that $H_2:=G_1[N_{G_1}(z)]$ satisfies $c<v(H_2)<2c+1$ and $\delta(H_2)\geqslant\lfloor c\rfloor$. Let $J$, $u$, $n$ and $G_0$ be given by Lemma~\ref{lem:red2}. By Lemma~\ref{lem:red2}, $n\in\{n_0+1,\dots,N\}$ and $\Phi(n)\leqslant\binom n2$, and so $(V(n),D(n))$ is certified for $t$ by hypothesis. Since $G_0$ is a non-null graph, $D(n)\leqslant\delta(G_0)\leqslant v(G_0)-1\leqslant V(n)-1$. So, by the definition of a certified pair, there is an integer $D'\leqslant D(n)$ and a sequence with charged costs that is admissible for $t$ with respect to $V(n)$ and $R':=2D'-V(n)+1$. By Lemma~\ref{lem:red2}, $v(G_0)\leqslant V(n)$ and
\[
2\delta(G_0)-v(G_0)\ \geqslant\ 2D(n)-V(n)\ \geqslant\ 2D'-V(n)\ =\ R'-1 .
\]
By Lemma~\ref{lem:core} with $m=t-2$, the graph $G_0$ has a $K_{t-2}$ minor. By Lemma~\ref{lem:red2}, $H_2$ has a $K_{t-1}$ minor, and by Lemma~\ref{lem:red1}, $G$ has a $K_t$ minor.
\end{proof}

\begin{lemma}\label{lem:small}
For every $t\in\{4,\dots,400\}$, every graph $G$ with $d(G)\geqslant4.5\,t\sqrt{\log t}$ has a $K_t$ minor.
\end{lemma}

\begin{proof}
Fix $t$. The script \texttt{exact\_recheck.py}, described in Appendix~\ref{sec:verification}, computes an integer $f\geqslant t-2$ and a function $\Phi:\{t-2,\dots,2f\}\to\mathbb Z$. Let $c:=f$, so that $\lfloor c\rfloor=f$ and $N=2f$, and let $n_0:=t-2\leqslant\lfloor c\rfloor$. The script checks that $\Phi$ is a density function for $c$ and $n_0$, and that $(V(n),D(n))$ is certified for $t$ for every $n\in\{t-1,\dots,2f\}$ with $\Phi(n)\leqslant\binom n2$. The script also checks that $f\leqslant4.5\,t\sqrt{\log t}$. So $d(G)\geqslant4.5\,t\sqrt{\log t}\geqslant f=c$, and Lemma~\ref{lem:assembly} shows that $G$ has a $K_t$ minor.
\end{proof}

For large $t$, we apply Lemma~\ref{lem:red2} with a density function given by a formula. Let $c>0$ be a real number. \defn{Slope data for $c$} consist of an integer $\nu\geqslant1$, positive integers $R_1,\dots,R_\nu$ and integers $Y_1<\dots<Y_{\nu-1}$ such that the numbers $r_j:=R_j+\frac12$ satisfy $r_1<2c+1$, and also $r_1<Y_1$ and $Y_{\nu-1}<2c+1$ if $\nu\geqslant2$. For such data, let $Y_0:=r_1$ and $Y_\nu:=2c+1$, let $v_0:=r_1$ and
\[
v_j\ :=\ v_{j-1}+r_j\log\frac{Y_j}{Y_{j-1}}\qquad\text{for }j\in\{1,\dots,\nu\} ,
\]
and let
\[
V_j\ :=\ \Bigl\lceil v_j+\frac23\Bigr\rceil-1\qquad\text{and}\qquad\tilde R_j\ :=\ \Bigl\lfloor r_j\Bigl(1-\frac1{2Y'_{j-1}}\Bigr)\Bigr\rfloor\qquad\text{for }j\in\{1,\dots,\nu\} ,
\]
where $Y'_0:=R_1+1$ and $Y'_j:=Y_j$ for $j\in\{1,\dots,\nu-1\}$. If $Y'_{j-1}\geqslant r_j$, then $r_j(1-\frac1{2Y'_{j-1}})\geqslant r_j-\frac12$, and so $\tilde R_j\geqslant R_j$, since $r_j-\frac12=R_j$ is an integer.

\begin{lemma}\label{lem:slopes}
For every real number $c>0$ and all slope data for $c$ with $v_\nu\leqslant\lfloor c\rfloor$, every graph $H_2$ with $c<v(H_2)<2c+1$ and $\delta(H_2)\geqslant\lfloor c\rfloor$ has a minor $J$ with a vertex $u$ such that the graph $G_0:=J[N_J(u)]$ is a non-null graph and satisfies $v(G_0)\leqslant V_j$ and $2\delta(G_0)-v(G_0)\geqslant\tilde R_j-1$ for some $j\in\{1,\dots,\nu\}$, and such that, for every integer $s\geqslant1$, if $G_0$ has a $K_s$ minor, then $H_2$ has a $K_{s+1}$ minor.
\end{lemma}

\begin{proof}
For real $y\geqslant r_1$, let $\omega(y):=v_{j-1}+r_j\log(y/Y_{j-1})$, where $j:=\nu$ if $y\geqslant Y_{\nu-1}$, and otherwise $j$ is the index with $Y_{j-1}\leqslant y<Y_j$. Since $\omega(Y_j)=v_j$ for $j\in\{1,\dots,\nu-1\}$, the function $\omega$ is continuous. Moreover, $\omega$ is increasing, and $\omega(r_1)=r_1$ and $\omega(2c+1)=v_\nu$. Let $\varphi(y):=y\,\omega(y)/2$, which is increasing. Let $n_0:=R_1+1$ and $\Phi(n):=\lceil\varphi(n)\rceil$ for $n\in\{n_0,\dots,N\}$. Since $r_1<2c+1$, we have $R_1<r_1<\omega(2c+1)=v_\nu\leqslant\lfloor c\rfloor$. So $n_0\leqslant\lfloor c\rfloor$.

We check that $\Phi$ is a density function for $c$ and $n_0$. Since $\varphi$ is increasing, $\Phi$ is non-decreasing. Since $n_0-1=R_1<r_1\leqslant\omega(n_0)$, we have $\Phi(n_0)\geqslant n_0\,\omega(n_0)/2>\binom{n_0}2$. For $n\in\{\lfloor c\rfloor+1,\dots,N\}$, we have $n<2c+1$, and so $\omega(n)<v_\nu\leqslant\lfloor c\rfloor$. Hence $\varphi(n)\leqslant n\lfloor c\rfloor/2$, which gives~\eqref{eq:density}.

Let $J$, $u$, $n$ and $G_0$ be given by Lemma~\ref{lem:red2}. Then $n\geqslant n_0+1\geqslant3$. Since $n\leqslant N<Y_\nu$, there is a least $j\in\{1,\dots,\nu\}$ with $n\leqslant Y_j$. If $j=1$, then $n-1\geqslant n_0=Y'_0>r_1=Y_0$. If $j\geqslant2$, then $n>Y_{j-1}$, and so $n-1\geqslant Y_{j-1}=Y'_{j-1}$, since both are integers. In both cases, $n-1\geqslant Y'_{j-1}$ and $[n-1,n]\subseteq[Y_{j-1},Y_j]$. So $\omega(n)-\omega(n-1)=r_j\log\frac n{n-1}$ and $\omega(n)\leqslant\omega(Y_j)=v_j$.

Since $V(n)\leqslant2\Phi(n)/n$, $\Phi(n)\geqslant\varphi(n)$ and $\Phi(n-1)<\varphi(n-1)+1$,
\begin{align*}
2D(n)-V(n)\ &\geqslant\ 2\Phi(n)\,\frac{n-1}n-2\Phi(n-1)\ >\ 2\varphi(n)\,\frac{n-1}n-2\varphi(n-1)-2\\
&=\ (n-1)\bigl(\omega(n)-\omega(n-1)\bigr)-2\ =\ r_j(n-1)\log\frac n{n-1}-2 .
\end{align*}
Since $\log(1+y)\geqslant y-y^2/2$ for $y\geqslant0$ and $n-1\geqslant Y'_{j-1}$, the right-hand side is at least $r_j(1-\frac1{2Y'_{j-1}})-2$. Since $2D(n)-V(n)$ is an integer, $2D(n)-V(n)\geqslant\tilde R_j-1$. Also $V(n)\leqslant2\Phi(n)/n<(2\varphi(n)+2)/n=\omega(n)+2/n\leqslant v_j+\frac23$. Since $V(n)$ is an integer, $V(n)\leqslant V_j$. By Lemma~\ref{lem:red2}, $v(G_0)\leqslant V(n)\leqslant V_j$ and $2\delta(G_0)-v(G_0)\geqslant2D(n)-V(n)\geqslant\tilde R_j-1$. The final claim is part of Lemma~\ref{lem:red2}.
\end{proof}

The next lemma builds a sequence $(k_i)$ for the core lemma from a small amount of data, so that the hypotheses of the core lemma reduce to finitely many inequalities. Let $V$ and $R\geqslant1$ be integers with $h:=\frac{V-R-1}2>0$, and let
\[
\pi_{\rm lo}\ :=\ \frac hV\qquad\text{and}\qquad\pi_{\rm hi}\ :=\ \frac h{V-R} .
\]
Then $0<\pi_{\rm lo}<\pi_{\rm hi}<\frac12$. \defn{Threshold data for $V$ and $R$} consist of integers $P\geqslant1$ and $2\leqslant\kappa_1<\dots<\kappa_P$, and real numbers $\pi_{\rm lo}=b_0<b_1<\dots<b_P=\pi_{\rm hi}$. For $\pi_{\rm lo}\leqslant\pi\leqslant\pi_{\rm hi}$, let $p(\pi)$ be the least index $p\in\{1,\dots,P\}$ with $\pi\leqslant b_p$, and let $\mu(\pi):=\kappa_{p(\pi)}+1+\lfloor\kappa_{p(\pi)}\pi\rfloor$. The \defn{threshold sequence} for these data is defined as follows. Let $S_0:=0$. For $i=1,2,\dots$, as long as $S_{i-1}\leqslant R$, let
\[
\pi_i\ :=\ \frac h{V-S_{i-1}},\qquad k_i\ :=\ \kappa_{p(\pi_i)}\qquad\text{and}\qquad S_i\ :=\ S_{i-1}+\mu(\pi_i).
\]
Since $0\leqslant S_{i-1}\leqslant R$, we have $\pi_{\rm lo}\leqslant\pi_i\leqslant\pi_{\rm hi}$. Since $k_i\geqslant2$, we have $\mu(\pi_i)=\ell_{k_i}(\pi_i)$. So $S_i$ and $\pi_i$ are the numbers of Section~\ref{sec:core} for the sequence $(k_i)$ with canonical costs. Position $i$ is \defn{in phase $p(\pi_i)$}. Since $\pi_i$ is increasing in $i$, the phase of position $i$ is non-decreasing in $i$.

For $p\in\{1,\dots,P\}$, let $M_p:=\kappa_p+1+\lfloor\kappa_pb_{p-1}\rfloor$, and for $1\leqslant q\leqslant p\leqslant P$, let $\eta_{qp}:=(\kappa_q-1)\kappa_p\geqslant2$. For $p\in\{1,\dots,P\}$, let
\[
A_p\ :=\ b_p^{\kappa_p}\Biggl[\,\sum_{q=1}^{p}\Bigl(1+\frac1{\kappa_q}\Bigr)^{\kappa_p}\frac{h\bigl(b_q^{\eta_{qp}-1}-b_{q-1}^{\eta_{qp}-1}\bigr)}{(\eta_{qp}-1)M_q}\ +\ \sum_{q=1}^{p-1}\Bigl(1+\frac1{\kappa_q}\Bigr)^{\kappa_p}b_q^{\eta_{qp}}\Biggr] ,
\]
and let
\[
\mathcal N\ :=\ h\int_{\pi_{\rm lo}}^{\pi_{\rm hi}}\frac{d\pi}{\pi^2\mu(\pi)} .
\]

\begin{lemma}\label{lem:threshold}
For all integers $V$ and $R\geqslant1$ with $V-R\geqslant2$, every integer $m\geqslant1$ and all threshold data for $V$ and $R$ such that $A_p<\frac1{\kappa_p+1}$ for $p\in\{1,\dots,P\}$ and $\mathcal N>m$, every non-null graph $G_0$ with $v(G_0)\leqslant V$ and $2\delta(G_0)-v(G_0)\geqslant R-1$ has a $K_{m+1}$ minor.
\end{lemma}

\begin{proof}
We first show that positions $1,\dots,m$ are defined and $S_m\leqslant R$. Suppose, for a contradiction, that some $i\leqslant m$ is such that position $i$ is not defined or $S_i>R$, and let $i$ be the least such index. Then $S_{i-1}\leqslant R$, so position $i$ is defined and $S_i>R$. Every position $a$ satisfies $h/\pi_a=V-S_{a-1}$. So $h(\frac1{\pi_a}-\frac1{\pi_{a+1}})=S_a-S_{a-1}=\mu(\pi_a)$ for $a\in\{1,\dots,i-1\}$, and $h(\frac1{\pi_i}-\frac1{\pi_{\rm hi}})=R-S_{i-1}<\mu(\pi_i)$. The function $\mu$ is non-decreasing on $[\pi_{\rm lo},\pi_{\rm hi}]$, since $p(\pi)$ and $\lfloor\kappa\pi\rfloor$ are non-decreasing in $\pi$ and $\kappa_1<\dots<\kappa_P$. So, for $\pi_{\rm lo}\leqslant x\leqslant y\leqslant\pi_{\rm hi}$,
\[
h\int_x^y\frac{d\pi}{\pi^2\mu(\pi)}\ \leqslant\ \frac h{\mu(x)}\int_x^y\frac{d\pi}{\pi^2}\ =\ \frac{h}{\mu(x)}\Bigl(\frac1x-\frac1y\Bigr) .
\]
Apply this inequality with $(x,y)=(\pi_a,\pi_{a+1})$ for $a\in\{1,\dots,i-1\}$, and with $(x,y)=(\pi_i,\pi_{\rm hi})$. Since $\pi_1=\pi_{\rm lo}$, summing the resulting inequalities gives $\mathcal N<i\leqslant m$, a contradiction.

We now check~\eqref{eq:core} at a position $i\leqslant m$. Let $p$ be the phase of $i$, and let $k:=\kappa_p$. Then $\pi_i\leqslant b_p$. Let $a<i$, let $q\leqslant p$ be the phase of $a$, and let $\eta:=\eta_{qp}$. Then $w_a^k=(1+\frac1{\kappa_q})^k\pi_a^{\eta}$. Suppose that $q=p$, or that $a$ is not the last position of phase $q$. Then position $a+1$ is in phase $q$, and so $b_{q-1}\leqslant\pi_a\leqslant\pi_{a+1}\leqslant b_q$. Write $x:=\pi_a$ and $y:=\pi_{a+1}$. Since $x\geqslant b_{q-1}$, the charged cost $\mu(x)=h(\frac1x-\frac1y)$ of position $a$ is at least $M_q$. Since $\eta\geqslant2$ and $0<x\leqslant y$,
\[
y^{\eta-1}-x^{\eta-1}\ =\ (y-x)\sum_{j=0}^{\eta-2}x^jy^{\eta-2-j}\ \geqslant\ (\eta-1)\,x^{\eta-2}(y-x) .
\]
Hence
\[
M_q\,x^{\eta}\ \leqslant\ h\,\frac{y-x}{xy}\,x^{\eta}\ \leqslant\ h\,x^{\eta-2}(y-x)\ \leqslant\ \frac{h\bigl(y^{\eta-1}-x^{\eta-1}\bigr)}{\eta-1} .
\]
The intervals $[\pi_a,\pi_{a+1}]$ for these positions $a$ have disjoint interiors and lie in $[b_{q-1},b_q]$. So the sum of $\pi_a^{\eta}$ over these positions $a$ is at most $h(b_q^{\eta-1}-b_{q-1}^{\eta-1})/((\eta-1)M_q)$. If $q<p$, then the last position $a$ of phase $q$ satisfies $\pi_a^{\eta}\leqslant b_q^{\eta}$. Multiplying these bounds by $(1+\frac1{\kappa_q})^k$, summing over $q\leqslant p$ and using $\pi_i^k\leqslant b_p^k$ gives $\sum_{a<i}(w_a\pi_i)^k\leqslant A_p$. By hypothesis, $A_p<\frac1{k+1}$, which gives~\eqref{eq:core} at $i$.

Therefore the sequence $k_1,\dots,k_m$ with canonical costs satisfies the hypotheses of Lemma~\ref{lem:core} with $m+1$ in place of $m$, since the sum of the canonical costs is $S_m\leqslant R<V$. So Lemma~\ref{lem:core} shows that $G_0$ has a $K_{m+1}$ minor.
\end{proof}

We now combine Lemmas~\ref{lem:slopes} and~\ref{lem:threshold} to handle the range $t\in\{401,\dots,10^{35}\}$. Lemma~\ref{lem:middle} below is proved by checking finitely many cells, each of which covers an interval of values of $t$. A \defn{cell} consists of integers $t_a\leqslant t_b$ and $f\geqslant1$, slope data for $f$, an integer $V_j^+\geqslant V_j$ for $j\in\{1,\dots,\nu\}$, and threshold data for $V_j^+$ and $\tilde R_j$ for $j\in\{1,\dots,\nu\}$. A cell is \defn{verified} if $f\leqslant4.5\,t_a\sqrt{\log t_a}$, $v_\nu\leqslant f$, and, for $j\in\{1,\dots,\nu\}$, the threshold data for $V_j^+$ and $\tilde R_j$ satisfy $A_p<\frac1{\kappa_p+1}$ for every $p$ and $\mathcal N>t_b-3$.

\begin{lemma}\label{lem:middle}
For every $t\in\{401,\dots,10^{35}\}$, every graph $G$ with $d(G)\geqslant4.5\,t\sqrt{\log t}$ has a $K_t$ minor.
\end{lemma}

\begin{proof}
The script \texttt{tail\_cells.py}, described in Appendix~\ref{sec:verification}, checks that each cell in the file \texttt{tail\_cells\_45.txt} is verified, and that every $t\in\{401,\dots,10^{35}\}$ satisfies $t\in\{t_a,\dots,t_b\}$ for some cell. Fix $t\in\{401,\dots,10^{35}\}$ and a cell with $t\in\{t_a,\dots,t_b\}$. Since $x\mapsto x\sqrt{\log x}$ is increasing for $x\geqslant1$, we have $d(G)\geqslant4.5\,t_a\sqrt{\log t_a}\geqslant f$. By Lemma~\ref{lem:red1} with $c:=f$, the graph $G$ has a minor $G_1$ with a vertex $z$ such that $H_2:=G_1[N_{G_1}(z)]$ satisfies $f<v(H_2)<2f+1$ and $\delta(H_2)\geqslant f$. By Lemma~\ref{lem:slopes} with $c:=f$, the graph $H_2$ has a minor $J$ with a vertex $u$ such that $G_0:=J[N_J(u)]$ is a non-null graph with $v(G_0)\leqslant V_j\leqslant V_j^+$ and $2\delta(G_0)-v(G_0)\geqslant\tilde R_j-1$ for some $j$. Since the cell is verified, the threshold data for $V_j^+$ and $\tilde R_j$ satisfy $\mathcal N>t_b-3\geqslant t-3$. So Lemma~\ref{lem:threshold} with $V:=V_j^+$, $R:=\tilde R_j$ and $m:=t-3$ shows that $G_0$ has a $K_{t-2}$ minor. By Lemmas~\ref{lem:slopes} and~\ref{lem:red1}, $G$ has a $K_t$ minor.
\end{proof}

For $\log t\geqslant80$, an analytic argument replaces the cells. The argument uses Lemma~\ref{lem:slopes} with three slopes, and Lemma~\ref{lem:threshold} with the canonical threshold data defined below. Let
\[
(\rho_1,\rho_2,\rho_3)\ :=\ (0.305,\,0.351,\,0.389)\qquad\text{and}\qquad(\sigma_1,\sigma_2)\ :=\ (0.485,\,0.903) .
\]
For an integer $t$ with $\lambda:=\log t\geqslant80$, let $c:=4.5\,t\sqrt\lambda$, let $R_j:=\lceil\rho_jc\rceil$ for $j\in\{1,2,3\}$, and let $Y_j:=\lceil\sigma_jc\rceil$ for $j\in\{1,2\}$. By~(H1) below, these numbers $R_j$ and $Y_j$ are slope data for $c$ with $\nu=3$. Let $r_j$, $v_j$, $V_j$ and $\tilde R_j$ be as defined before Lemma~\ref{lem:slopes}.

For integers $V$ and $R\geqslant1$ with $h:=\frac{V-R-1}2\geqslant1$, let $\pi_{\rm lo}$ and $\pi_{\rm hi}$ be as defined before Lemma~\ref{lem:threshold}, let $\Lambda_h:=\log(2\mathrm eh)$, and, for every integer $k\geqslant2$, let
\[
\Pi_k\ :=\ \Bigl(\frac{k^2-k-1}{2\mathrm eh}\Bigr)^{1/(k^2-1)},\qquad\text{so that}\qquad(k^2-1)\log\frac1{\Pi_k}\ =\ \Lambda_h-\log(k^2-k-1).
\]
Since $\Pi_k\to1$ as $k\to\infty$, there is a least integer $k_{\rm f}\geqslant2$ with $\Pi_{k_{\rm f}}>\pi_{\rm lo}$, and a least integer $k_{\rm l}\geqslant2$ with $\Pi_{k_{\rm l}}\geqslant\pi_{\rm hi}$. Since $\pi_{\rm hi}>\pi_{\rm lo}$, we have $k_{\rm f}\leqslant k_{\rm l}$. Let $P:=k_{\rm l}-k_{\rm f}+1$, let $\kappa_p:=k_{\rm f}+p-1$ for $p\in\{1,\dots,P\}$, let $b_p:=\Pi_{\kappa_p}$ for $p\in\{1,\dots,P-1\}$, and let $b_0:=\pi_{\rm lo}$ and $b_P:=\pi_{\rm hi}$. We call $P$, $(\kappa_p)$ and $(b_p)$ the \defn{canonical threshold data} for $V$ and $R$.

We show that the canonical threshold data are threshold data for $V$ and $R$. Let $k\in\{2,\dots,k_{\rm l}-2\}$. Then $\Pi_{k+1}<\pi_{\rm hi}<1$ by the choice of $k_{\rm l}$, so $a_{k+1}:=\Lambda_h-\log((k+1)^2-(k+1)-1)>0$. Also $a_k:=\Lambda_h-\log(k^2-k-1)>a_{k+1}$, since $k^2-k-1<(k+1)^2-(k+1)-1$. So
\[
\log\frac1{\Pi_k}\ =\ \frac{a_k}{k^2-1}\ >\ \frac{a_{k+1}}{k^2-1}\ >\ \frac{a_{k+1}}{(k+1)^2-1}\ =\ \log\frac1{\Pi_{k+1}} ,
\]
that is, $\Pi_k<\Pi_{k+1}$. Hence $b_1<\dots<b_{P-1}$. If $P\geqslant2$, then $b_1=\Pi_{k_{\rm f}}>\pi_{\rm lo}$, and $b_{P-1}=\Pi_{k_{\rm l}-1}<\pi_{\rm hi}$ by the choice of $k_{\rm l}$, since $k_{\rm l}-1\geqslant k_{\rm f}\geqslant2$. If $P=1$, then $b_0=\pi_{\rm lo}<\pi_{\rm hi}=b_1$. Also $b_P=\pi_{\rm hi}\leqslant\Pi_{\kappa_P}$, and so $b_p\leqslant\Pi_{\kappa_p}$ for every $p\geqslant1$.

Let $u:=1/c$. Then $0<u\leqslant u_{\max}:=1/(4.5\,\mathrm e^{80}\sqrt{80})<10^{-36}$. For $j\in\{1,2,3\}$, let $(V,R):=(V_j,R_j)$ and $h:=\frac{V-R-1}2$. The script \texttt{tail\_hand.py}, described in Appendix~\ref{sec:verification}, computes, for all $u\in(0,u_{\max}]$ simultaneously, a lower bound $\alpha_j$ for $h/c$ and bounds $p_j^-\leqslant\pi_{\rm lo}\leqslant p_j^+$. To four decimal places, $\alpha_1,\alpha_2,\alpha_3$ are $0.0707$, $0.1568$, $0.2925$, and $p_j^-$ and $p_j^+$ are $0.1584$, $0.2359$, $0.3003$ for $j=1,2,3$. Let $k^*_j:=6+j$ and $\beta_j:=\frac12(1+\frac1{k^*_j})(1+\frac1{k^*_j+1})$. For $0<\pi<1$, let
\[
\gamma(\pi)\ :=\ (1+\pi)\sqrt{\frac1{\log(1/\pi)}+\frac1{80}}+\frac{2+\pi}{\sqrt{80}}\qquad\text{and}\qquad\Gamma_j\ :=\ \int_{p_j^+}^{0.4999}\frac{d\pi}{\pi^2\gamma(\pi)} .
\]
The script \texttt{tail\_hand.py} checks the following facts, which are used in the proof of Lemma~\ref{lem:large}, where $j$ ranges over $\{1,2,3\}$ in~(H3)--(H6).
\begin{enumerate}[label=\textup{(H\arabic*)},nosep]
\item $r_1<Y_1<Y_2<2c$, $Y_1\geqslant r_2$ and $Y_2\geqslant r_3$.
\item $v_3\leqslant c-2$.
\item $4.5\,\alpha_j\sqrt{80}\geqslant1$ and $9\mathrm e\alpha_j>1$.
\item $4.5\,\alpha_j\Gamma_j>\sqrt{1+(\frac12\log80+\log(9\mathrm e\alpha_j))/80}$.
\item $((k^*_j-1)^2-1)\log(1/p_j^-)+\log((k^*_j-1)^2-k^*_j)\leqslant80$.
\item $2^{-k^*_j-1}\beta_j/((1-\beta_j)\,p_j^-)+(k^*_j+2)^2\,2^{-(k^*_j+1)(k^*_j-1/2)}\leqslant\frac12$.
\end{enumerate}
\begin{lemma}\label{lem:large}
For every integer $t$ with $\log t\geqslant80$, every graph $G$ with $d(G)\geqslant4.5\,t\sqrt{\log t}$ has a $K_t$ minor.
\end{lemma}

\begin{proof}
Let $\lambda:=\log t$ and $c:=4.5\,t\sqrt\lambda$. By Lemma~\ref{lem:red1}, $G$ has a minor $G_1$ with a vertex $z$ such that $H_2:=G_1[N_{G_1}(z)]$ satisfies $c<v(H_2)<2c+1$ and $\delta(H_2)\geqslant\lfloor c\rfloor$. By~(H1), the numbers $R_j$ and $Y_j$ are slope data for $c$, and $\tilde R_j\geqslant R_j$ for every $j$. By~(H2), $v_3\leqslant c-2<\lfloor c\rfloor$. So Lemma~\ref{lem:slopes} shows that $H_2$ has a minor $J$ with a vertex $x$ such that $G_0:=J[N_J(x)]$ is a non-null graph with $v(G_0)\leqslant V_j$ and $2\delta(G_0)-v(G_0)\geqslant R_j-1$ for some $j$. Fix this $j$, and let $V:=V_j$, $R:=R_j$, $h:=\frac{V-R-1}2$ and $m:=t-3$. By Lemma~\ref{lem:threshold}, it suffices to show that the canonical threshold data for $V$ and $R$ are defined and satisfy $A_p<\frac1{\kappa_p+1}$ for every $p$ and $\mathcal N>m$. Indeed, Lemma~\ref{lem:threshold} then shows that $G_0$ has a $K_{t-2}$ minor, and so $G$ has a $K_t$ minor by Lemmas~\ref{lem:slopes} and~\ref{lem:red1}.

By the choice of $\alpha_j$ and by~(H3), $h\geqslant\alpha_jc=4.5\,\alpha_jt\sqrt\lambda\geqslant t$. So $\Lambda_h\geqslant\log h\geqslant\lambda\geqslant80$. Since also $R\geqslant1$, the canonical threshold data are defined. By the choice of $p_j^-$ and $p_j^+$, we have $p_j^-\leqslant\pi_{\rm lo}\leqslant p_j^+$. Also $\pi_{\rm hi}=\frac12-\frac1{2(V-R)}\geqslant0.4999$, since $V-R=2h+1>5000$.

We first show that $k_{\rm f}\geqslant k^*_j$. Let $k\in\{2,\dots,k^*_j-1\}$. Then $\Pi_k\leqslant\pi_{\rm lo}$ if and only if $(k^2-1)\log(1/\pi_{\rm lo})+\log(k^2-k-1)\leqslant\Lambda_h$. The left-hand side is non-decreasing in $k$, and $\log(1/\pi_{\rm lo})\leqslant\log(1/p_j^-)$. Since $(k^*_j-1)^2-(k^*_j-1)-1=(k^*_j-1)^2-k^*_j$ and $\Lambda_h\geqslant80$, fact~(H5) gives this inequality for $k=k^*_j-1$, and hence $\Pi_k\leqslant\pi_{\rm lo}$ for every such $k$. So $k_{\rm f}\geqslant k^*_j\geqslant7$.

We next show that $A_p<\frac1{\kappa_p+1}$ for every $p$. Let $p\in\{1,\dots,P\}$ and $k:=\kappa_p$. Every breakpoint is at most $\pi_{\rm hi}<\frac12$. In the term $q=p$ of the first sum in $A_p$, drop the subtracted term, and use $\eta_{pp}=k^2-k$, $b_p\leqslant\Pi_k$, $M_p\geqslant k+1$ and $(1+\frac1k)^k<\mathrm e$. Then this term, multiplied by $b_p^k$, is at most
\[
\Bigl(1+\frac1k\Bigr)^k\frac{h\,\Pi_k^{k^2-1}}{(k^2-k-1)(k+1)}\ =\ \Bigl(1+\frac1k\Bigr)^k\frac1{2\mathrm e(k+1)}\ <\ \frac1{2(k+1)} .
\]
Now let $q<p$, let $k':=\kappa_q$, let $d:=k-k'\geqslant1$ and let $\eta:=\eta_{qp}=(k'-1)k$. Then $b_q=\Pi_{k'}$, so $\pi_{\rm lo}<b_q<1$ and $h\,b_q^{k'^2-1}=(k'^2-k'-1)/(2\mathrm e)$. Since $\eta-k'^2=(k'-1)(d-1)-1\geqslant-1$, it follows that $h\,b_q^{\eta-1}\leqslant h\,b_q^{k'^2-2}<(k'^2-k'-1)/(2\mathrm e\pi_{\rm lo})$. Since $\eta-1\geqslant k'^2-2\geqslant k'^2-k'-1$, $M_q\geqslant k'+1$ and $b_p<\frac12$, the term $q$ of the first sum in $A_p$, multiplied by $(k+1)b_p^k$, is at most
\[
\frac{(k+1)(1+\frac1{k'})^k}{2^{k+1}\mathrm e\,\pi_{\rm lo}(k'+1)}\ \leqslant\ \frac{2^{-k'-1}}{\pi_{\rm lo}}\biggl[\frac{(1+\frac1{k'})(1+\frac1{k'+1})}2\biggr]^d\ \leqslant\ \frac{2^{-k^*_j-1}}{p_j^-}\,\beta_j^d .
\]
Here the first inequality uses $(1+\frac1{k'})^k<\mathrm e(1+\frac1{k'})^d$ and $\frac{k+1}{k'+1}\leqslant(1+\frac1{k'+1})^d$, and the second uses $k'\geqslant k_{\rm f}\geqslant k^*_j$. Distinct $q<p$ give distinct $d\geqslant1$, so the terms $q<p$ of the first sum contribute at most $2^{-k^*_j-1}\beta_j/((1-\beta_j)p_j^-)$ in total. Since $(1+\frac1{k'})^k\leqslant(\frac43)^k\leqslant2^{k/2}$, the term $q$ of the second sum in $A_p$, multiplied by $(k+1)b_p^k$, is
\[
(k+1)\Bigl(1+\frac1{k'}\Bigr)^kb_p^kb_q^\eta\ \leqslant\ (k+1)\,2^{k/2}\,2^{-k}\,2^{-(k'-1)k}\ \leqslant\ (k+1)\,2^{-k(k^*_j-1/2)} .
\]
There are $p-1<k$ such terms. Since $q<p$ exists only if $k\geqslant k_{\rm f}+1\geqslant k^*_j+1$, and $(k+1)^22^{-k(k^*_j-1/2)}$ is decreasing in $k$, the terms of the second sum contribute at most $(k^*_j+2)^22^{-(k^*_j+1)(k^*_j-1/2)}$ in total. Adding these three bounds and using~(H6) gives $(k+1)A_p<\frac12+\frac12=1$.

Finally, we show that $\mathcal N>m$. Let $\pi_{\rm lo}\leqslant\pi\leqslant\pi_{\rm hi}$, let $p:=p(\pi)$ and let $k:=\kappa_p$. We claim that $k\leqslant1+\sqrt{1+\Lambda_h/\log(1/\pi)}$. The claim is clear if $k=2$. If $k\geqslant3$, then $\pi\geqslant\Pi_{k-1}$. Indeed, if $p\geqslant2$, then $\pi>b_{p-1}=\Pi_{k-1}$, and if $p=1$, then $k=k_{\rm f}$ and $\Pi_{k_{\rm f}-1}\leqslant\pi_{\rm lo}\leqslant\pi$ by the choice of $k_{\rm f}$. Since $(k-1)^2-(k-1)-1\geqslant1$, it follows that $((k-1)^2-1)\log(1/\pi)\leqslant((k-1)^2-1)\log(1/\Pi_{k-1})\leqslant\Lambda_h$, which proves the claim. So, using $\Lambda_h\geqslant80$,
\[
\mu(\pi)\ \leqslant\ k(1+\pi)+1\ \leqslant\ (2+\pi)+(1+\pi)\sqrt{1+\frac{\Lambda_h}{\log(1/\pi)}}\ \leqslant\ \sqrt{\Lambda_h}\,\gamma(\pi).
\]
Since $\pi_{\rm lo}\leqslant p_j^+$ and $0.4999\leqslant\pi_{\rm hi}$, it follows that $\mathcal N\geqslant h\,\Gamma_j/\sqrt{\Lambda_h}$. The function $y\mapsto y/\sqrt{\log(2\mathrm ey)}$ is increasing for $y\geqslant1$, and $h\geqslant\alpha_jc\geqslant1$. Since $\log(2\mathrm e\alpha_jc)=\lambda+\frac12\log\lambda+\log(9\mathrm e\alpha_j)$,
\[
\mathcal N\ \geqslant\ \frac{\alpha_jc\,\Gamma_j}{\sqrt{\log(2\mathrm e\alpha_jc)}}\ =\ t\cdot\frac{4.5\,\alpha_j\Gamma_j}{\sqrt{1+(\frac12\log\lambda+\log(9\mathrm e\alpha_j))/\lambda}} .
\]
By~(H3), $\log(9\mathrm e\alpha_j)>0$, and so the function $\lambda\mapsto(\frac12\log\lambda+\log(9\mathrm e\alpha_j))/\lambda$ is decreasing for $\lambda\geqslant\mathrm e$. Since $\lambda\geqslant80$, fact~(H4) gives $\mathcal N>t>m$.
\end{proof}

Lemmas~\ref{lem:small}, \ref{lem:middle} and~\ref{lem:large} together prove Theorem~\ref{thm:Kt}.

\begin{proof}[Proof of Theorem~\ref{thm:Kt}]
Let $G$ be a graph with $d(G)\geqslant4.5\,t\sqrt{\log t}$. If $t=2$, then $d(G)>0$, so $G$ has an edge, which gives a $K_2$ minor. If $t=3$, then $d(G)>1$, since $4.5\cdot3\sqrt{\log3}>1$. So $G$ is not a forest and thus contains a cycle, which gives a $K_3$ minor. If $t\in\{4,\dots,400\}$, then Lemma~\ref{lem:small} applies. If $t\in\{401,\dots,10^{35}\}$, then Lemma~\ref{lem:middle} applies. If $t>10^{35}$, then $\log t>80$, since $\mathrm e^{80}<5.6\cdot10^{34}$, and Lemma~\ref{lem:large} applies.
\end{proof}

%% Generated by make_bib.py from ../refs.bib. Do not edit by hand.

%=====================================================================
\appendix
\section{Computer verification}\label{sec:verification}
%=====================================================================

This appendix describes the scripts \texttt{exact\_recheck.py}, \texttt{tail\_cells.py} and \texttt{tail\_hand.py}, which are used in the proofs of Lemmas~\ref{lem:small}, \ref{lem:middle} and~\ref{lem:large}, respectively.

\paragraph{\boldmath The range $t\in\{4,\dots,400\}$.} We now describe the script \texttt{exact\_recheck.py}. For every integer $V\geqslant t-2$, a floating-point search, on which the proof does not rely, gives an integer $\overline D(V)$, and $\overline D(V):=V$ for $V<t-2$. The script sets $\Phi(t-2):=\binom{t-2}2+1$. For $n\geqslant t-1$, the script defines $\Phi(n)$ to be the least integer $X\geqslant\Phi(n-1)$ such that
\[
X>\binom n2,\qquad\text{or}\qquad D\geqslant V,\qquad\text{or}\qquad D\geqslant\overline D(V),
\]
where $V:=\lfloor2X/n\rfloor$ and $D:=X-\Phi(n-1)$. So $\Phi$ is non-decreasing and $\Phi(t-2)>\binom{t-2}2$. The script takes $f$ to be the least integer $f\geqslant t-2$ for which~\eqref{eq:density} with $\lfloor c\rfloor=f$ and $N=2f$ holds at every $n\in\{f+1,\dots,2f\}$, and the script checks these inequalities in integer arithmetic. So $\Phi$ is a density function for $f$ and $t-2$. The script also checks that $f\leqslant4.5\,t\sqrt{\log t}$, in interval arithmetic.

Now let $n\in\{t-1,\dots,2f\}$ with $\Phi(n)\leqslant\binom n2$. If $D(n)\geqslant V(n)$, then $(V(n),D(n))$ is certified by definition. Otherwise, let $V:=V(n)$. Since $\Phi(n)\leqslant\binom n2$ and $D(n)<V$, the choice of $X=\Phi(n)$ gives $D(n)\geqslant\overline D(V)$. So $V\geqslant t-2$, since $\overline D(V')=V'$ for $V'<t-2$. For every such $V$, the script checks that a sequence $k_1,\dots,k_{t-3}$ with canonical costs is admissible for $t$ with respect to $V$ and $R:=2\overline D(V)-V+1$. So $(V(n),D(n))$ is certified for $t$, with $D':=\overline D(V)$. Admissibility is checked in exact integer arithmetic, with each term $w_a^{k_i}$ in~\eqref{eq:core} rounded up to an integer multiple of $2^{-512}$.

Every check passes, and the largest value of $f/(t\sqrt{\log t})$ is $4.46302$, at $t=115$ and $f=1118$.

\paragraph{\boldmath The range $t\in\{401,\dots,10^{35}\}$.} We now describe the script \texttt{tail\_cells.py}, which imports \texttt{tail\_core.py}. In its \texttt{check} mode, the script reads the cells from the file \texttt{tail\_cells\_45.txt} and checks the following, independently of how the cells were found.
\begin{itemize}[nosep]
\item The first cell has $t_a=401$, the last cell has $t_b\geqslant10^{35}$, and every other cell has $t_a$ equal to $t_b+1$ for the previous cell.
\item The inequality $f\leqslant4.5\,t_a\sqrt{\log t_a}$, in interval arithmetic with $200$ bits of precision, using the \texttt{mpmath} library.
\item The conditions $f\geqslant1$, $R_j\geqslant1$ for every $j$, $r_1<2f+1$, $Y_1<\dots<Y_{\nu-1}$, $r_1<Y_1$ and $Y_{\nu-1}<2f+1$, and the inequality $v_\nu\leqslant f$, where each $v_j$ is evaluated in interval arithmetic. The integer $V_j^+$ is $\lceil x\rceil-1$, where $x$ is the upper end of the interval computed for $v_j+\frac23$, and so $V_j^+\geqslant V_j$. The integers $\tilde R_j$ are computed exactly.
\item For each $j$, that $V_j^+-\tilde R_j\geqslant2$, that the values $\kappa_p$ are integers with $2\leqslant\kappa_1<\dots<\kappa_P$, and that the breakpoints $b_p$ are rational and satisfy $\pi_{\rm lo}<b_1<\dots<b_P=\pi_{\rm hi}$. Since $R_j\geqslant1$ and $Y'_{j-1}\geqslant2$, also $\tilde R_j\geqslant1$. The script checks the inequalities $(\kappa_p+1)A_p<1$ for every $p$ and $\mathcal N\geqslant t_b-2$ in exact rational arithmetic, using the Python module \texttt{fractions}. Since $\mu$ is a step function, $\mathcal N$ is a finite sum of terms $h(\frac1x-\frac1y)/\mu_0$, where $\mu_0$ is the constant value of $\mu$ on an interval $(x,y)$.
\end{itemize}
Every check passes. The file contains $728$ cells, and the last cell has $t_b=10^{35}$. Each cell has $\nu\in\{8,\dots,120\}$, and every value $\kappa_p$ is at most $11$. The cells were found by floating-point searches, in the \texttt{gen} mode of the script, with the breakpoints $b_p$ for $p<P$ rounded down to integer multiples of $2^{-50}$, and the proof does not rely on these searches.

\paragraph{\boldmath The range $\log t\geqslant80$.} We now describe the script \texttt{tail\_hand.py}, which checks the numerical facts used in the proof of Lemma~\ref{lem:large}. The script evaluates all quantities in interval arithmetic with $120$ bits of precision, using the \texttt{mpmath} library, with $u$ ranging over $[0,u_{\max}]$. Since $R_j\in[\rho_jc,\rho_jc+1)$, $Y_j\in[\sigma_jc,\sigma_jc+1)$ and $V_j\in[v_j-\frac13,v_j+\frac23)$, the script encloses $r_j/c$, $R_j/c$, $Y_j/c$ and $V_j/c$ in intervals that contain $[\rho_j+\frac u2,\rho_j+\frac{3u}2]$, $[\rho_j,\rho_j+u]$, $[\sigma_j,\sigma_j+u]$ and $[v_j/c-\frac u3,v_j/c+\frac{2u}3]$, respectively. The script computes $v_j/c$ from $r_j/c$, $Y_j/c$ and $(2c+1)/c=2+u$, and then $h/c=(V_j/c-R_j/c-u)/2$ and $\pi_{\rm lo}=(h/c)/(V_j/c)$. The number $\alpha_j$ is the lower end of the interval computed for $h/c$, and $p_j^-$ and $p_j^+$ are the ends of the interval computed for $\pi_{\rm lo}$. The script bounds $\Gamma_j$ from below by a right-endpoint Riemann sum with $4000$ steps, which is valid since the integrand is decreasing in $\pi$. The two sides of~(H4) are $1.0231$ and $1.0170$ for $j=1$, $1.0297$ and $1.0219$ for $j=2$, and $1.0421$ and $1.0257$ for $j=3$. The left-hand sides of~(H6) are $0.0444$, $0.0138$ and $0.0051$ for $j=1,2,3$, respectively. The script takes $k^*_j$ to be the largest integer for which~(H5) holds, and this gives $k^*_j=6+j$. Every check passes.

The ancillary files of this paper include the file \texttt{tail\_cells\_45.txt} and the scripts \texttt{exact\_recheck.py}, \texttt{tail\_cells.py}, \texttt{tail\_core.py} and \texttt{tail\_hand.py}. We ran the scripts with version $1.4.1$ of the \texttt{mpmath} library.

\paragraph{Formal verification.} Theorem~\ref{thm:Kt} has been formally verified in the proof assistant Lean~4 with the Mathlib library. The Lean files are among the ancillary files on arXiv. The main Lean statement is the following.
\begin{verbatim}
theorem KtMinor45.kt_minor_45 :
  ∀ (t : ℕ), 2 ≤ t → ∀ (V : Type u) [Fintype V]
  (G : SimpleGraph V) [DecidableRel G.Adj], 1 ≤ Fintype.card V →
  (4.5 : ℝ) * t * √(Real.log t) * Fintype.card V ≤ G.edgeFinset.card →
  HasKMinor G t
\end{verbatim}
Here \texttt{HasKMinor G t} says that $G$ has a model of $K_t$, as defined in Section~\ref{sec:intro}. The Lean kernel checks every computation in the proof. These computations check the conditions in the proofs of Lemmas~\ref{lem:small}, \ref{lem:middle} and~\ref{lem:large}, with the following differences from the scripts. First, in the proof of Lemma~\ref{lem:small}, each sequence is given by runs of equal values with equal charged costs. The inequality $L_i\geqslant\ell_{k_i}(\pi_i)$ and condition~\eqref{eq:core} are checked only at the last position of each run, and the part of the sum in~\eqref{eq:core} over each run is bounded by an integral, which is evaluated exactly in rational arithmetic. Since admissibility with respect to $V$ and $R$ implies admissibility with respect to $V'\leqslant V$ and $R'\geqslant R$ with $R'<V'$, only $46\,109$ certificates are needed for all $t\in\{4,\dots,400\}$. The density function is computed from these certificates by a greedy rule like that of \texttt{exact\_recheck.py}. This density function differs from the one computed by \texttt{exact\_recheck.py}, and it gives a slightly larger integer $f$ for some values of $t$. However, $f\leqslant4.5\,t\sqrt{\log t}$ still holds for every $t\in\{4,\dots,400\}$, and the largest value of $f/(t\sqrt{\log t})$ is $4.49995$, at $t=267$. Second, logarithms and square roots are bounded by rational numbers, instead of by interval arithmetic. In the proof of Lemma~\ref{lem:middle}, the integers $V_j^+$ are read from the file \texttt{tail\_cells\_45.txt} and checked, instead of being recomputed. Third, in the proof of Lemma~\ref{lem:large}, the numbers $\alpha_j$, $p_j^-$ and $p_j^+$ are replaced by rational numbers with eight decimal places. These numbers are shown to satisfy $h\geqslant\alpha_jc$ and $p_j^-\leqslant\pi_{\rm lo}\leqslant p_j^+$ using only $c\geqslant2^{80}$. The integral $\Gamma_j$ is bounded below by a sum over $100$ intervals. Facts (H1)--(H5) are checked for these rational numbers, and~(H6) is replaced by a similar inequality that comes from a slightly different bound on the second sum in $A_p$.

\end{document}